# Estimation errors of the Sharpe ratio for long-memory stochastic volatility models

Hwai-Chung Ho[1]

*Academia Sinica*

**Abstract:** The Sharpe ratio, which is defined as the ratio of the excess expected return of an investment to its standard deviation, has been widely cited in the financial literature by researchers and practitioners. However, very little attention has been paid to the statistical properties of the estimation of the ratio. Lo (2002) derived the $\sqrt{n}$-normality of the ratio's estimation errors for returns which are iid or stationary with serial correlations, and pointed out that to make inference on the accuracy of the estimation, the serial correlation among the returns needs to be taken into account. In the present paper a class of time series models for returns is introduced to demonstrate that there exists a factor other than the serial correlation of the returns that dominates the asymptotic behavior of the Sharpe ratio statistics. The model under consideration is a linear process whose innovation sequence has summable coefficients and contains a latent volatility component which is long-memory. It is proved that the estimation errors of the ratio are asymptotically normal with a convergence rate slower than $\sqrt{n}$ and that the estimation deviation of the expected return makes no contribution to the limiting distribution.

## 1. Introduction

An interesting phenomenon observed in many financial time series is that strong evidence of persistent correlation exists in some nonlinear transformation of returns, such as square, logarithm of square, and absolute value, whereas the return series itself behaves almost like white noise. This so-called clustering volatility property has a profound implication. The traditional linear processes such as ARMA models and the mixing conditions of various types that have been widely used to account for the weak-dependence or short-memory properties of stationary processes (see, e.g., [1]) are found inadequate to model the dependence structure of the return process. A great deal of research works have been devoted to looking for proper models that entail the stylized fact mentioned above. The ARCH model proposed by Engle [6] and its various extensions are attempts that have been proved very successful. Recently, models other than ARCH family have been seen to provide better fitting for data with clustering volatility. For instance, Lobato and Savin [11] examine the S&P 500 index series for the period of July 1962 to December 1994 and report that the squared daily returns exhibit the genuine long-memory effect which ARCH process cannot produce (see also [5]). Based on Lobato and Savin's finding, Breidt, Crato and Lima [2] suggest the following long-memory stochastic volatility model (LMSV):

$$(1.1) \qquad r_t = v_t \varepsilon_t, \quad v_t = \delta \exp(x_t),$$

[1]Institute of Statistical Science, Academia Sinica, Taipei 115, Taiwan, e-mail: hcho@stat.sinica.edu.tw

where $\delta > 0$, and $\{x_t\}$ is a Gaussian process which exhibits long memory and is independent of the iid sequence $\{\varepsilon_t\}$ with mean zero and variance one. The short-memory version of model (1.1) has been discussed, for example, by Taylor [14], Melino and Turnbull [12] and Harvey et al. [8]. The precise definition of short- or long-memory process is given as follows. A linaer process defined as

$$(1.2) \qquad x_t = \sum_{i=0}^{\infty} a_i z_{t-i},$$

where the $z_i$ are iid random variables (Gaussian or non-Gaussian) having mean 0 and variance one, is called short-memory if the coefficients $a_i$ are summable or long-memory if $a_i \sim C i^{-\beta}$ with $\beta$ in $(1/2, 1)$; "$g_n \sim h_n$" signifies $\lim_{n \to \infty} g_n/h_n = 1$. The long-memory process just defined is sometimes also referred to as a fractional differencing (or $I(d)$) process with the memory parameter $d = 1 - \beta$ [3]. It can be seen that the LMSV model described in (1.1) and (1.2) exhibits the desirable property that $\{r_t\}$ is white noise and $\{r_t^2\}$ is long-memory. Because of this characteristic property, one needs to be cautious in making statistical inference for the LMSV model if the statistics of interest involve nonlinear transformations. The purpose of this paper is to point out a circumstance under which the estimation statistics based on the LMSV model behave distinctly different from traditional stationary sequences of weak dependence such as the ARMA model with iid innovations.

We use the example of the Sharpe ratio to demonstrate that for the LMSV model the estimation statistics have entirely different asymptotic properties from those of the case where the volatilily is short-memory. Discussions of this and a more general model are presented in Sections 2 and 3, respectively. The paper's main result is formulated in a theorem stated in Section 3 and its proof is given in Section 4.

## 2. LMSV models: the simple case

The Sharpe raio, which is defined as the ratio of the excess expected return of an investment to its standard deviation, is originally motivated by the mean-variance analysis and the Sharpe-Lintner Captial Asset Pricing Model (Campbell, Lo and MacKinlay [4]) and has become a popular index used to evaluate investment performance and for risk management. Both the expected return and the standard deviation are generally unknown and need to be estimated. Although the ratio is one of the most commonly cited statistics in financial analysis by researchers and practitioners as well, not much attention has been paid to its statistical properties until the work of Lo [10]. Lo [10] points out that to gauge the accuracy of the estimates of the ratio, it is important to take into account the dependence of the returns for it may result significant difference of the limiting variance between iid and non-iid (dependent) returns. For both of the two cases the standard $\sqrt{n}$ central limit theorem is assumed to hold for the ratio's estimates. The LMSV time series is a stationary martingale difference sequence bearing strong dependence in the latent component of volatility. The partial sums of the sequence itself and of the sequence after a certain transformation is applied may have entirely different asymptotic behaviors. Below we show that for the LMSV model, the Sharpe ratio statistic is asymptotically normal but converges to the true ratio at a rate slower than $\sqrt{n}$. Furthermore, while the ratio's statistics involve the estimates of the expected return and the standard deviation, it turns out that only the estimation errors of the latter contribute to the limit distribution as opposed to the case of short-memory volatility where neither of the two estimates is asymptotically negligible.

Let the returns $\{r_t\}$ be model as in (1.1) and (1.2) with long-memory $x_t = \sum_{i=0}^{\infty} a_i z_{t-i}$, where $z_i$ are iid random variables having mean 0 and variance 1 and the coefficients $a_i$ are such that $a_i \sim C \cdot i^{-\beta}$ with $\beta$ in $(1/2, 1)$. Denote by $\sigma^2 = Er_t^2$. For the observed returns $\{r_1, \ldots, r_n\}$, we define

$$\hat{\mu} = n^{-1} \sum_{t=1}^{n} r_t, \quad \hat{\sigma}^2 = n^{-1} \sum_{t=1}^{n} (r_t - \hat{\mu})^2,$$

and the Sharpe ratio statistics

$$\hat{SR} = \frac{\hat{\mu} - r_f}{\hat{\sigma}},$$

where $r_f$ is a fixed risk-free interest rate assumed to be positive. Using the $\delta$-method, we have

$$\hat{SR} - SR = \frac{\hat{\mu}}{\sigma} + \frac{r_f(\hat{\sigma}^2 - \sigma^2)}{2\sigma^3} + O_p((\hat{\sigma}^2 - \sigma^2)^2).$$

Also write

$$\hat{\sigma}^2 - \sigma^2 = n^{-1} \sum_{t=1}^{n} v_t^2(\varepsilon_t^2 - 1) + n^{-1} \sum_{t=1}^{n} (v_t^2 - \sigma^2) - \hat{\mu}^2.$$

To derive the asymptotic distribution, we first compute the variance of $\hat{\sigma}^2 - \sigma^2$. Note that

(2.1) $$\text{var}(n^{-1} \sum_{t=1}^{n} v_t^2(\varepsilon_t^2 - 1)) = O(n^{-1}) \quad \text{and} \quad \text{var}(\hat{\mu}) = O(n^{-1}),$$

since both $\{v_t^2(\varepsilon_t^2 - 1)\}$ and $\{v_t \varepsilon_t\}$ are sequences of martingale differences. For $\sum_{t=1}^{n}(v_t^2 - \sigma^2)$, we use the results obtained by Ho and Hsing [9]. Let $F(\cdot)$ be the common distribution function of the $x_t$. Denote by

$$K_{\infty}(y) = e^{2y} \int e^{2x} dF(x).$$

Then by Theorem 3.1 and Corollary 3.3 of Ho and Hsing [9],

(2.2) $$n^{\beta - 3/2} \{\sum_{t=1}^{n}(v_t^2 - \sigma^2)\} = \delta^2 K_{\infty}^{(1)}(0) n^{\beta - 3/2} \{\sum_{t=1}^{n} x_t\} + o_p(1)$$

$$\xrightarrow{d} 2\sigma^2 \cdot N(0, \xi^2)$$

with

$$\xi^2 = C^2 \frac{\int_0^{\infty}(x^2 + x)^{-\beta} dx}{2(1-\beta)(3/2+\beta)} \cdot \int_{-\infty}^{1} \{\int_0^1 [(v-u)^+]^{-\beta} dv\} du.$$

Combining (2.1) and (2.2) gives

(2.3) $$n^{\beta - 1/2}(\hat{SR} - SR) = \frac{r_f}{2\sigma^3} n^{\beta - 1/2}(\hat{\sigma}^2 - \sigma^2) + o_p(1)$$

$$\xrightarrow{d} r_f \sigma^{-1} \cdot N(0, \xi^2),$$

If $x_t$ is short-memory in the sense as specified before that

$$x_t = \sum_{i=0}^{\infty} a_i z_{t-i} \quad \text{with} \quad \sum_{i=1}^{\infty} |a_i| < \infty,$$

then the usual $\sqrt{n}$ central limit theorem will hold for $\sqrt{n}(\hat{SR} - SR)$. The proof of this will be covered in the next subsection as a special case of a more general model.

## 3. Linear processes of LMSV models

We now focus on the linear process with its innovations being a LMSV sequence. Specifically, define

$$(3.1) \qquad y_t = \sum_{j=0}^{\infty} b_j r_{t-j} \quad \text{with} \quad \sum_j |b_j| < \infty,$$

where $r_t$ is modeled in (1.1) and (1.2) with $\delta = 1$. Denote by $\sigma_y^2$ the variance of the $y_t$. The Sharpe ratio now is $SR = r_f/\sigma_y$ and its corresponding estimator is

$$(3.2) \qquad \hat{SR}_y = \frac{W_n - r_f}{\hat{\sigma}_y},$$

where

$$W_n = n^{-1} \sum_{t=1}^{n} y_t, \quad \hat{\sigma}_y = (n^{-1} \sum_{t=1}^{n} (y_t - W_n)^2)^{1/2}.$$

From now on we assume that there is a positive constant $K$ such that for any $\eta > 0$

$$(3.3) \qquad Ee^{\eta x_1} \leq e^{K\eta^2}.$$

As can be seen later in the proof we only need a sufficiently large constant $K$. Using a stronger condition here is merely for the ease of presentation.

**Theorem.** *For the model defined in (3.1), assume condition (3.3) holds.*

(i) *Suppose $x_t$ is short-memory, that is, $\sum_{i=0}^{\infty} |a_i| < \infty$. Assume $E\varepsilon_1^3 = 0$, then*

$$(3.4) \qquad \sqrt{n}(\hat{SR} - SR) \xrightarrow{d} N(0, \xi_1^2)$$

*for some constant $\xi_1$.*

(ii) *If $x_t$ is long-memory with the coefficients satisfying that $a_i \sim Ci^{-\beta}$ for $\beta \in (1/2, 1)$, then*

$$(3.5) \qquad n^{\beta - 3/2}(\hat{SR} - SR) \xrightarrow{d} 2 \int e^{2x} dF(x) N(0, \xi_2^2)$$

*for some constant $\xi_2$.*

The limiting variances, $\xi_1^2$ and $\xi_2^2$, given in (3.4) and (3.5) above will be derived in the proof of the theorem. Both $\xi_1^2$ and $\xi_2^2$ depend on the linear filter $\{b_j\}$ and some parameters of the laten process $\{x_t\}$. It is a very challanging problem to estimate the two quantities. For part (ii) of the Theorem, if the distribution function $F(\cdot)$ of $x_t$ is known, then one can use the sampling window method proposed in [7] and [13] to consistently estimate $\xi_1^2$ and $\xi_2^2$. As for the short-memory case of part (i) of the Theorem, no existing results in the literature cover this case unless a certain kind of weak dependence is assumed. With only the summability condition on $\{a_j\}$ one needs to develop some new theory to support the use of the resampling scheme mentioned above.

*Proof of Theorem.* (i) Define

$$x_{t,m} = \sum_{i=0}^{m-1} a_i \varepsilon_{t-i}, \quad \tilde{x}_{t,m} = \sum_{i=m}^{\infty} a_i \varepsilon_{t-i}, \quad r_{t,m} = e^{x_{t,m}} \varepsilon_t, \quad y_{t,m} = \sum_{j=0}^{m-1} b_j r_{t-j,m},$$

$$W_{n,m} = n^{-1} \sum_{t=1}^{n} y_{t,m}.$$

Since $y_{t,m}$'s are $2m$-dependent, as $n \to \infty$,

(3.6)
$$\sqrt{n} W_{n,m} \xrightarrow{d} N(0, \lambda_m^2),$$

where

$$\lambda_m^2 = \lim_{n \to \infty} n^{-1} \mathrm{var}(\sum_{t=1}^{n} y_{t,m})$$

$$= \delta^2 E e^{2x_{1,m}} (\sum_{j=0}^{m} b_j^2 + 2 \sum_{k=1}^{\infty} \sum_{j=0}^{m} b_j b_{j+k}).$$

Write

$$\sqrt{n}(W_n - W_{n,m}) = n^{-1/2} \sum_{t=1}^{n} (y_t - y_{t,m})$$

$$= n^{-1/2} \sum_{t=1}^{n} \sum_{j=0}^{m-1} b_j (r_{t-j} - r_{t-j,m}) + n^{-1/2} \sum_{t=1}^{m} \sum_{j=m}^{\infty} b_j r_{t-j}$$

$$\equiv C_{n,m} + D_{n,m}.$$

Then

$$E C_{n,m}^2 = \delta^2 E e^{2x_{1,m}} (e^{\tilde{x}_{1,m}} - 1) \sum_{k=-n-1}^{n-1} (1 - \frac{|k|}{n})(\sum_{j=0}^{m-1} b_j b_{j+k}).$$

By using the elementary inequality $|e^x - 1| \leq e|x|$, $|x| \leq 1$, and the Chebyshev inequality, we have

$$E(e^{\tilde{x}_{1,m}} - 1)^2 = E(e^{\tilde{x}_{1,m}} - 1)^2 I\{\tilde{x}_{1,m} \leq 1\} + E(e^{\tilde{x}_{1,m}} - 1)^2 I\{\tilde{x}_{1,m} > 1\}$$

$$\leq e(E\tilde{x}_{1,m}^2) + (E(e^{\tilde{x}_{1,m}} - 1)^4)^{1/2} (E\tilde{x}_{1,m}^2)^{1/2}.$$

Because, by assumption (3.3), $Ee^{4x_{1,m}}$ is bounded in $m$, we have

(3.7)
$$E(e^{\tilde{x}_{1,m}} - 1)^2 \to 0 \quad \text{as} \quad m \to \infty.$$

This and $\sum_{j}^{\infty} |b_j| < \infty$ jointly imply

(3.8)
$$\lim_{m \to \infty} \lim_{n \to \infty} EC_{n,m}^2 = 0.$$

Similarly,

$$\lim_{m\to\infty}\lim_{n\to\infty} ED_{n,m}^2 = 0. \tag{3.9}$$

From (3.6), (3.8) and (3.9) it follows that

$$\sqrt{n}W_n \to N(0, \lambda^2), \tag{3.10}$$

where

$$\lambda^2 = \lim_{m\to\infty} \lambda_m^2 = \sigma^2 \Big(\sum_{j=0}^{\infty} b_j^2 + \sum_{k=1}^{\infty}\sum_{j=0}^{\infty} b_j b_{j+k}\Big).$$

We now derive the limiting distribution for $\sqrt{n}(\hat{\sigma}_y^2 - \sigma_y^2)$. Write

$$\sqrt{n}(\hat{\sigma}_y^2 - \sigma_y^2) = \delta^2 n^{-1/2} \sum_{t=1}^{n}\sum_{j=0}^{\infty} b_j^2 e^{2x_{t-j}}(\varepsilon_{t-j}^2 - 1) + \delta^2 n^{-1/2} \sum_{t=1}^{n}\sum_{j=0}^{\infty}(e^{2x_{t-j}} - \sigma_y^2)$$

$$+ n^{-1/2} \sum_{t=1}^{n}\sum_{i\neq j} b_i b_j r_{t-i} r_{t-j}$$

$$\equiv V_{n,1} + V_{n,2} + V_{n,3}. \tag{3.11}$$

By the same $m$- truncation argument as used in proving (3.8) one can show that $V_{n,1}$, $V_{n,2}$ and $V_{n,3}$ are asymptotically normal and independent, that is, as $n \to \infty$,

$$V_{n,1} + +V_{n,2} + V_{n,3} \xrightarrow{d} N(0, g^2), \tag{3.12}$$

where $g^2$ is the sum of the limiting variances of $V_{n,1}$, $V_{n,2}$ and $V_{n,3}$. Because $x_t$ may be non-Gaussian, the analytic form of the covariance function of $\{e^{2x_t} - \sigma_y^2\}$ and consequently of the limiting variance of $V_{n,2}$, which equals to

$$\delta^2 \lim_{m\to\infty} \limsup_{n\to\infty} n^{-1} \text{var}\Big(\sum_{t=1}^{n}(e^{2x_{t,m}} - \sigma_y^2)\Big),$$

is not available. However, the exact formulas of limiting variances for $V_{n,1}$ and $V_{n,3}$ can be found as follows.

$$\lim_{n\to\infty} \text{var}(V_{n,1}) = \delta^4 [Ee^{4x_1}][E(\varepsilon_1^2 - 1)^2]\Big[\sum_{j=0}^{\infty} b_j^4 + 2\sum_{k=1}^{\infty}\sum_{j=0}^{\infty} b_j^2 b_{j+k}^2\Big],$$

$$\lim_{n\to\infty} \text{var}(V_{n,3}) = \delta^4 [Ee^{2x_1}]^2 \Big[\sum_{i\neq j} b_i^2 b_j^2 + \sum_{k=1}^{\infty}\sum_{i\neq j} b_i b_{i+k} b_j b_{j+k}\Big].$$

Note that the assumption $E\varepsilon_1^3 = 0$ is used to prove that $V_{n,2}$ is asymptotically independent with $V_{n,1}$ and $V_{n,3}$. The limit results of (3.10) and (3.12) imply

$$\sqrt{n}(\hat{SR} - SR) \xrightarrow{d} N(0, \xi_1^2)$$

with $\xi_1^2 = \lambda^2 + r_f^2 (4\sigma_y^6)^{-1} g^2$. Hence (2.4) holds.
(ii) Because $\{r_t\}$ is a sequence of martingale differences, we have

$$\text{var}(W_n) = O(1/n). \tag{3.13}$$

Similarly, for $V_{n,1}$ and $V_{n,3}$ defined in (3.11),

(3.14) $$\text{var}(V_{n,1}) = O(1), \quad \text{var}(V_{n,3}) = O(1).$$

To compute the variance of $V_{n,2}$, define $y'_t = \sum_{j=0}^{\infty} b_j x_{t-j}$. Then $y'_t$ can be rewritten as

$$y'_t = \sum_{j=0}^{\infty} z_{t-j} B_j,$$

where

$$B_j = \sum_{i=0}^{j} b_i a_{j-i}.$$

As $j \to \infty$, since $a_j \sim Cj^{-\beta}$,

$$B_j \sim C_1 j^{-\beta},$$

where $C_1 = C(\sum_{i=0}^{\infty} b_i)$. Then, as $k \to \infty$,

$$\sum_{j=0}^{\infty} B_j B_{j+k} \sim C_1^2 \int x^{-\beta}(1+x)^{-\beta} dx \cdot k^{-2\beta+1},$$

implying that

$$E y'_t y'_{t+k} = \sum_{j=0}^{\infty} B_j B_{j+k} \sim C_1^2 \int x^{-\beta}(1+x)^{-\beta} dx \cdot k^{-2\beta+1}.$$

In other words, $\{y'_t\}$ is also a linear long-memory process having the same memory parameter as that of $x_t$. Therefore, similar to (2.2),

(3.15) $$n^{\beta-3/2} \sum_{t=1}^{n} y'_t \xrightarrow{d} N(0, \xi_2^2)$$

with

$$\xi_2^2 = \frac{C_1^2 \int_0^{\infty}(x^2+x)^{-\beta} dx}{2(1-\beta)(3/2+\beta)} \cdot \int_{-\infty}^{1} \{\int_0^1 [(v-u)^+]^{-\beta} dv\} du.$$

As noted before in (2.2) that

$$n^{\beta-3/2}\{\sum_{t=1}^{n}(e^{2x_t} - \sigma^2)\} = 2\int e^{2x} dF(x) \cdot (n^{\beta-3/2}\{\sum_{t=1}^{n} x_t\}) + o_p(1).$$

From this and (3.15), we have, as $n \to \infty$,

$$n^{\beta-3/2}\{\sum_{t=1}^{n}\sum_{j=0}^{\infty} b_j(e^{2x_{t-j}} - \sigma_y^2)\} = 2\int e^{2x} dF(x)(n^{\beta-3/2}\sum_{t=1}^{n} y'_t) + o_p(1)$$

(3.16) $$\xrightarrow{d} 2\int e^{2x} dF(x) \cdot N(0, \xi_2^2).$$

Summarizing (3.13), (3.14) and (3.16) gives

$$n^{\beta-3/2}(\hat{SR} - SR) \xrightarrow{d} 2\int e^{2x} dF(x) N(0, \xi_2^2).$$

The proof is completed. □